\documentclass[12pt,a4paper]{article}
\usepackage{amssymb}
\usepackage{amsmath}

\newtheorem{theorem}{Theorem}

\newtheorem{corollary}[theorem]{Corollary}

\newtheorem{definition}[theorem]{Definition}

\newtheorem{lemma}[theorem]{Lemma}

\newenvironment{proof}[1][Proof]{\textbf{#1.} }{\ \rule{0.5em}{0.5em}}

\newcommand\E{\mathbf{E}}

\newcommand{\calP}{{\cal P}}
\newcommand{\ep}{\varepsilon}
\newcommand{\dN}{{{\mathbb N}}}
\newcommand{\dR}{{{\mathbb R}}}

\newcommand{\conv}{{\rm{conv}}}

\newcommand{\calM}{\mathcal{M}}

\usepackage{graphicx}

\newcommand{\numbercellong}[2]
{
\begin{picture}(40,20)(0,0)
\put(0,0){\framebox(40,20)} \put(20,10){\makebox(0,0){#1}}
\put(33,15){\makebox(0,0){#2}}
\end{picture}
}

\newcounter{figurecounter}
\begin{document}

\title{Bounded Computational Capacity Equilibrium%
\thanks{This work was conducted while the second author was visiting Universidad de Valencia.
The first author thanks both the Spanish Ministry of Science
and Technology and the European Feder Founds for financial support under
project SEJ2007-66581 and Generalitat Valenciana (PROMETEO/2009/068).
The second author thanks the Departamento de An\'alisis Econ\'omico at Universidad de Valencia
for the hospitality during his visit.
The authors thank Elchanan Ben Porath, Ehud Kalai and Ehud Lehrer for their suggestions.
The work of Solan was partially supported by ISF grant 212/09.}}
\author{\normalsize Pen\'elope Hern\'andez%
\thanks{ERI-CES and Departamento de An\'alisis Econ\'omico, Universidad de Valencia.
Campus de Los Naranjos s/n, 46022 Valencia, Spain.}
and Eilon Solan%
\thanks{Department of Statistics and Operations Research, School of Mathematical Sciences, Tel Aviv University,
 Tel Aviv 69978, Israel. eilons@post.tau.ac.il}}
\maketitle

\begin{abstract}
We study repeated games played by players with bounded computational power,
where, in contrast to Abreu and Rubisntein \cite{AR88},
the memory is costly.
We prove a folk theorem: the limit set of equilibrium payoffs in mixed strategies,
as the cost of memory goes to 0, includes the set of feasible and individually rational payoffs.
This result stands in sharp contrast to \cite{AR88},
who proved that when memory is free,
the set of equilibrium payoffs in repeated games played by players with bounded computational power
is a strict subset of the set of feasible and individually rational payoffs.
Our result emphasizes the role of memory cost and of mixing when players have bounded computational power.
\end{abstract}

\thispagestyle{empty}


\section{Introduction}

In a seminal work, Simon \cite{S72}, \cite{S78} recognized the impact of bounded rationality in economic modelization
of individual agents and of organizations.
In the last few decades an expanding literature studied the implementation cost of strategies in strategic interactions
(see, e.g., Rubinstein \cite{R98}, Chatterjee and Sabourian \cite{CS08}).
One particular research question deals with
achieving a target outcome as a collusion, cooperation, or bargaining outcome by non-sophisticated agents
(see, e.g., Chatterjee and Sabourian \cite{CS00}, Sabourian \cite{S03}, Gale and Sabourian \cite{GS05},
and Maenner \cite{M08}).

One common way to model players with bounded rationality is by restricting them to
strategies that can be implemented by finite state machines, or automata.
The game theoretic literature on repeated games played by finite automata can be roughly divided into two categories.
On the one hand, an extensive literature (e.g., Kalai \cite{Kalai},
Ben Porath \cite{B93}, Piccione \cite{P92}, Piccione and Rubinstein \cite{PR93},
Neyman \cite{N85}, \cite{N97}, \cite{N98},   Neyman and Okada \cite{NO99}, \cite{NO00}, \cite{NOIJGT00}, Zemel \cite{Z89})
study games where the memory size of the two players is determined exogenously,
so that each player can deviate only to strategies with the given memory size.
On the other hand, Rubinstein \cite{R86},
Abreu and Rubinstein \cite{AR88} and Banks and Sundaram \cite{BS90} study games where the players have lexicographic
preferences: each player tries to maximize her payoff, and subject to that she tries to minimize her memory size.
Thus, it is assumed that memory is free, and a player would deviate to significantly more
complex strategy if that would increase her profit by one cent.

In practice, the level of complexity that players can use in their strategies is not known in advance,
either because players do not know each other computational power,
because players may increase their computational power if they realize that such an increase is beneficial,
or because players may decrease their computational power if the loss caused by this decrease
is compensated by the reduced expenses due to this decision.

In the present paper we take a more pragmatic point of view than the two approaches mentioned above,
and we study repeated games played by boundedly rational players,
when the computational power is costly.

As a motivating example,
consider employees' training for a new job.
The training period enables the employee to cope with situations that he may encounter in the future.
The longer the training period, the better prepared will be the employee,
thereby increasing the employer's profit.

Once the employee starts working, he follows the instructions that he learned,
and so we can model the employee as a finite state machine.
The training period dictates the size of the machine,
that is, the number of its states,
and the training itself determines how the machine behaves in various situations.
Because training is costly,
the employer will try to balance between the length of the training period and the gains from extended training.

When employees of different employers interact, say a salesperson and a buyer,
the evolution of the interaction is dictated by their training.
A salesperson, say, may interact with buyers of different firms,
who undertook different training programs,
and therefore follow different finite state machines.
Therefore he has some uncertainty regarding the finite state machine that the buyer will follow,
so that in fact he faces a mixed strategy.

The employers, who plan the training of their respective employees,
then face a game,
where each tries to teach her employees the techniques that best cope with the techniques
taught by the other employer.
As salespersons and buyers interact repeatedly,
the situation can be modelled as a repeated game played by finite state machines,
where the goal of each player is to maximize some combination of the long-run average payoff
and the cost of training.

To capture situations like the one in the example, we
assume for simplicity that the players have additive utility:
the utility of a player is the sum of her long-run average payoff
and the cost of her computational power.
Formally, for every positive real number $c$, we say that the vector $x \in \dR^2$
is a {\em $c$-Bounded Computational Capacity equilibrium} (hereafter, BCC for short) if
it is an equilibrium when the utility of each player is the difference between her long-run
average payoff and $c$ times the size of its finite state machine.

A payoff vector $x \in \dR^2$ is a {\em BCC equilibrium payoff}
if it is the limit, as $c$ goes to 0,
of payoffs that correspond to $c$-bounded computational capacity payoffs,
and the cost of the machines used along the sequence converges to 0.

Interestingly, the definition does not imply that the set of BCC equilibrium payoffs
is a subset, nor a super set, of the set of Nash equilibrium payoffs.

Our main result is a folk theorem:
in two player games,
every feasible and individually rational (w.r.t. the min-max value in pure strategies)
payoff vector is a BCC equilibrium payoff.

Our proof is constructive: we explicitly construct equilibrium strategies.
The equilibrium play is composed of three phases.
The first phase, that on the equilibrium path is played only once,
is a punishment phase;
in this phase each player plays a strategy that punishes the other player,
that is, an action that attains the min-max value in pure strategies of the opponent.
As in \cite{AR88}, it is crucial to have the punishment phase on the equilibrium path;
otherwise, players can use smaller machines, that cannot implement punishment and lower the cost of their machines.
However, if a machine cannot implement punishment,
there is nothing that will deter the other player from deviating.
The second phase, called the babbling phase,
is also played only once on the equilibrium path.
In this phase the players play a predetermined sequence of action pairs.
In the third phase, called the regular phase,
the players play repeatedly a predetermined periodic sequence of action pairs
that approximates the desired target payoff.
To implement this phase, the players re-use states that were used in the babbling phase.
In fact, the role of the babbling phase is to enable one to embed the regular phase within it,
and its structure is designed to simplify complexity calculations.
It is long enough to ensure that with only low probability a player can
correctly guess which of the states in the other player's machine are re-used.

One can describe the equilibrium path by imagining the following meeting between two strangers.
At first, the strangers exchange threats and vivid descriptions of what each one will do to the other
if the other does not behave as desired.
After they prove to each other that they can execute punishment,
they indulge in a long small-talk.
Finally, they go to business, and implement the desired outcome.

Our paper is closely related to Abreu and Rubinstein \cite{AR88},
where a characterization of the set of equilibrium payoffs is provided when the players have a lexicographic preference:
subject to maximizing her long-run average payoff each player wishes to minimize the
complexity of the finite state machine that
implements her strategy. The main result of \cite{AR88} is that
the set of equilibrium payoffs is the set of all feasible and individually rational payoffs
(relative to the min-max value in pure strategies)
that can be supported by coordinated play.
The main message of \cite{AR88} is that a folk theorem does not obtain:
the set of equilibrium payoffs may be strictly smaller than the set of feasible and individually rational payoffs.
Our result shows that two properties of the model of \cite{AR88} drive their result.
First, \cite{AR88} assumes that computational power is costless,
so that players will deviate to a prohibitively large automaton
to gain a cent.
This is in contrast to our model, where computational power is costly.
Second, \cite{AR88} restricts the players to pure strategies,
whereas we allow the players to use mixed strategies.

Abreu and Rubinstein \cite{AR88} point at a difficulty in using mixed strategies in games played by
players with bounded computational power: mixing is a complex operation,
and players with bounded computational power will prefer to use a pure strategy than a mixed strategy,
thereby saving the cost of mixing.
We argue that there are at least two interpretations of the model
where the use of mixed strategies is natural.
First, it may happen that the agent playing  the game
is limited, whereas the player who chooses the strategy for the agent does not have limits on her
computational power.
Thus, the complexity of computing the strategy played by the agent can be large,
and include mixing,
while the complexity of implementing this strategy should be low.
Second, a player may not know the identity of the agent whom her own agent is going to face,
and therefore she does not know the pure simple strategy which that agent is going to use.
Alternatively, the other agents who her agent is going to face may use different pure strategies.
Thus, the player may assume that the other player randomly chooses her simple strategy.
In our construction the role of  mixing is
to hide the strategy that each agent uses.
Whereas in \cite{AR88} the players use pure strategies to reduce their computational power,
which leads to a significantly smaller set of equilibrium payoffs,
mixing allows the players to use once again complex strategies,
and the folk theorem is restored.

The rest of the paper is organized as follows. Section 2 presents the model and the main result.
The construction of a mixed equilibrium strategy for both players in the particular case of the Prisoner's Dilemma is
presented in Section 3. In Section 4 we explain how the construction is adapted for general two-player games.

\section{The Model and the Main Result}

In this section we define the model, including the concepts of
automata, repeated games, and strategies implementable by an automaton;
we describe our solution concept of Bounded Computational Capacity equilibrium,
and we state the main result.

\subsection{Repeated Games}

A two-player {\em repeated game} is given by (1) two finite action
sets $A_1$ and $A_2$ for the two players, and (2) two payoff
functions $u_1 : A_1 \times A_2 \to \dR$ and $u_2 : A_1 \times A_2
\to \dR$ for the two players.

The game is played as follows. At every stage $t$, each player $i \in \{1,2\}$
chooses an action $a_i^t \in A_i$, and receives the stage payoff
$u_i(a_1^t,a_2^t)$. The goal of each player is to maximize its
long-run average payoff $\lim_{t \to
\infty}\frac{1}{t}\sum_{j=1}^t u_i(a_1^j,a_2^j)$,
where $\{(a_1^j,a_2^j), j \in \dN\}$ is the sequence of action pairs that were chosen by the players.%
\footnote{In general this limit need not exist. Our solution concept will
take care of this issue.}
A {\em pure strategy} of player $i$ is a function that assigns an action in $A_i$ to every finite history $h \in \cup_{t=0}^\infty (A_1 \times A_2)^t$.
A {\em mixed strategy} of player $i$ is a probability distribution over pure strategies.

\subsection{Automata}

A common way to model a decision maker with bounded computational capacity is as an automaton,
which is a finite state machine whose output depends on the current state,
and whose evolution depends on the current state and on its input
(see, e.g., Neyman \cite{N85} and Rubinstein \cite{R86}).
Formally, an {\em automaton} $P$ is given by (1) a finite state space $Q$,
(2) a finite set $I$ of inputs, (3) a finite set $O$ of outputs,
(4) an output function $f : Q \to O$, (5) a transition function
$g : Q \times I \to Q$, and (6) an initial state $q^* \in Q$.

Denote by $q^t$ the automaton's state at stage $t$.
The automaton starts in state $q^1 = q^*$, and at every stage $t$, as a
function of the current state $q^t$ and the current input $i^t$,
the output of the automaton $o^t = f(q^t)$ is determined, and the automaton moves to a
new state $q^{t+1} = g(q^t,i^t)$.

The {\em size} of an automaton $P$, denoted by $|P|$, is the
number of states in $Q$.
Below we will use strategies that can be implemented by automata;
in this case the size of the automaton measures the complexity of the strategy.

\subsection{Strategies Implemented by Automata}

Fix a player $i \in \{1,2\}$.
An automaton $P$ whose set of inputs is the set of actions of player $3-i$
and set of outputs is the set of actions of player $i$, that is, $I = A_{3-i}$ and $O = A_i$,
can implement a pure strategy of player $i$.
Indeed, at every stage $t$, the strategy
plays the action $f(q^t)$, and the new state of the automaton
$q^{t+1}=g(q^t,a_{3-i}^t)$ depends on its current state $q^t$ and on the action
$a_{3-i}^t$ that the other player played at stage $t$.
For $i=1,2$, we denote an automaton that implements a strategy of player $i$ by
$P_i$.
We denote by $\calP^m_i$ the set of all automata with $m$ states that
implement pure strategies of player $i$.

When the players use arbitrary strategies, the long-run average
payoff needs not exist. However, when both players use strategies
that can be implemented by automata, say $P_1$ and $P_2$
of sizes $p_1$ and $p_2$ respectively, the evolution of the
automata follows a Markov chain with $p_1 \times p_2$ states, and
therefore the long-run average payoff exists. We denote this
average payoff by $\gamma(P_1,P_2) \in \dR^2$.

A mixed automaton $M$ is a probability distribution over pure
automata%
\footnote{To emphasize the distinction between automata and mixed automata,
we call the former {\em pure automata}.}.
A mixed automaton corresponds to the situation in which
the automaton that is used is not known, and there is a belief
over which automaton is used.
A mixed automaton defines a mixed
strategy: at the outset of the game, a pure automaton is chosen
according to the probability distribution given by the mixed
automaton, and the strategy that the pure automaton defines is
executed.

We will use only mixed
automata whose support is pure automata of a given size $m$. Denote by
$\calM^m_i$ the set of all mixed automata whose support is
automata in $\calP^m_i$, and by $\calM_i = \cup_{m \in \dN} \calM_i^m$ the set of all mixed automata
whose support contains automata of the same size.
If $M_i \in \calM^m_i$, we say that $m$ is the size of the automaton $M_i$.
Thus, the size of a mixed automaton refers to the {\em size} of the pure automata in its support
(and not, for example, to the {\em number} of pure automata in its support).
If we interpret each pure automaton as an agent's type,
and a mixed automaton as the type's distribution in the population,
then the size of the mixed automaton measures the complexity of an individual agent,
and not the type diversity in the population.

When both players use mixed strategies that can be implemented by mixed automata,
the expected long-run average payoff exists;
it is the expectation of the long-run average payoff of the (pure) automata that the players play:
\[ \gamma(M_1,M_2) := \E_{M_1,M_2}[\gamma(P_1,P_2)]. \]

\subsection{Bounded Computational Capacity Equilibrium}

In the present section we study games where the utility function of each player takes into account
the complexity of the strategy that she uses.

\begin{definition}
Let $c > 0$.
A pair of mixed automata $(M_1,M_2)$ is a {\em $c$-BCC equilibrium},
if it is a Nash equilibrium for the utility functions $U^c_i(M_1,M_2) = \gamma_i(M_1,M_2) - c|M_i|$, $i\in \{1,2\}$.
\end{definition}

If the game has an equilibrium in pure strategies, then the pair of pure automata $(P_1,P_2)$,
both with size 1,
that repeatedly play the equilibrium actions of the two players,
is a $c$-BCC equilibrium, for every $c > 0$.

The min-max value of player $i$ in pure strategies in the one-shot game is
\[ v_i := \min_{a_{3-i} \in A_{3-i}} \max_{a_i \in A_i} u_i(a_i,a_{3-i}). \]
An action $a_{3-i}$ that attains the minimum is termed a {\em punishing} action of player $3-i$.

To get rid of the dependency of the constant $c$ we define the concept of a {\em BCC equilibrium payoff}.
A payoff vector $x$ is a {\em BCC equilibrium payoff} if it is the limit, as $c$ goes to 0,
of the payoff that corresponds to $c$-BCC equilibria.
\begin{definition}
A payoff vector $x = (x_1,x_2)$ is a {\em BCC equilibrium payoff} if for every $c > 0$ there is a
$c$-BCC equilibrium $(M_1(c),M_2(c))$ such that $\lim_{c \to 0} U^c(M_1(c),M_2(c)) = x$
and $\lim_{c \to 0} cM_i(c) = 0$.
\end{definition}

The condition $\lim_{c \to 0} U^c(M_1(c),M_2(c)) = x$ in the definition of a BCC-equilibrium payoff
ensures that $x$ can be supported as an equilibrium,
while the condition $\lim_{c \to 0} cM_i(c) = 0$ ensure that the cost of the automata that support
this equilibrium is negligible.
In particular, the limit of the long-run average payoffs also converges to $x$: $\lim_{c \to 0} \gamma(M_1(c),M_2(c)) = x$.

It follows from the discussion above that every pure equilibrium payoff is a BCC equilibrium payoff.
Using Abreu and Rubinstein's \cite{AR88} proof,
one can show
that any individually rational payoff (relative to the min-max value in pure strategies)
that can be generated by coordinated play is a BCC equilibrium payoff.
For the formal statement, assume w.l.o.g. that $|A_1| \leq |A_2|$.
\begin{theorem}[Abreu and Rubinstein, 1988]
\label{theorem ar}
Let $\sigma : A_1 \to A_2$ be a one-to-one function.
Then any payoff vector $x$ in the convex hull of $\{u(a_1,\sigma(a_1)), a_1 \in A_i\}$
that satisfies $x_i > v_i$ for $i=1,2$ is a BCC equilibrium payoff.
\end{theorem}

\subsection{The Main Result}

The set of feasible payoff vectors is
\[ F := \conv\{ u(a), a \in A_1 \times A_2\}. \]
The set of {\em strictly individually rational payoff vectors} (relative to the min-max value in pure strategies) is
\[ V := \left\{x = (x_1,x_2) \in \dR^2 \colon x_1 > v_1, x_2 > v_2\right\}. \]

Our main result is the following folk theorem,
that states that every feasible and strictly individually rational payoff vector is a BCC equilibrium payoff.
\begin{theorem}
\label{th1}
Every vector in $F \cap V$ is a BCC equilibrium payoff.
\end{theorem}

Observe that Theorem \ref{th1} is not a characterization of the set of BCC equilibrium payoffs,
because it does not rule out the possibility that a feasible payoff that is not individually rational
(relative to the min-max value in pure strategies) is a BCC equilibrium payoff.
That is, we do not know whether threats of punishments by a mixed strategy in the one-shot game
can be implemented in a BCC equilibrium.

Theorem \ref{th1} stands in sharp contrast to the main message of Abreu and Rubinstein  \cite{AR88},
where it is proved that lexicographic preferences,
which is equivalent to an infinitesimal cost function $c$, implies that in equilibrium players follow coordinated play,
so that the set of equilibrium payoffs is sometimes smaller than the set of feasible and individually rational payoffs.
Our study shows that the result of Abreu and Rubinstein  \cite{AR88}
hinges on two assumptions: (a) memory is costless,
and (b) the players use only pure automata.
Once we assume that memory is costly,
and that players may use mixed automata,
the set of equilibrium payoffs dramatically changes.

\subsection{Comments and Discussion}

\subsubsection{On the definition of BCC equilibria}

The definition of BCC equilibrium is analog to the definition of Nash equilibrium;
in both we ask whether a specific behavior (that is, a pair of strategies) is stable.
Thus, in a $c$-BCC equilibrium we assume that each player already has an automaton with which
she is going to play the game, and we ask whether playing this automaton is the best response given the
automaton that the other player is going to use.
As in the definition of Nash equilibrium,
we do not ask how the players arrived at these automata, and we do not restrict the sizes of these automata
(though the memory cost does bound the maximal size of automaton that the players will use).
In principle it may well be that some BCC equilibrium payoff can be supported only with prohibitively large automata,
which we would like to rule out.
That is, we may want to add the size of the automata that the players use to the definition itself.
In our construction (see the proof of Theorem \ref{th1}),
to support a $c$-BCC equilibrium payoff that is close to some target payoff $x$
we use two automata of similar sizes;
the size of the automaton is related to both $c$ and to the level of approximation to the target payoff:
as $c$ gets closer to 0, and as the $c$-BCC equilibrium payoff gets closer to $x$,
we use larger automata.

\subsubsection{BCC equilibria and Nash equilibria}

Theorem \ref{th1} states that every feasible and individually rational (w.r.t. the max-min value in pure strategies)
payoff vector is a BCC equilibrium payoff. This theorem does not rule out the possibility that there would be a
payoff vector that is {\em not} individually rational that would still be a BCC equilibrium;
that is, a BCC equilibrium payoff need not be a Nash equilibrium payoff.
The theorem also does not rule out the possibility that some payoff vector that is
individually rational w.r.t. the max-min value in {\em mixed} strategies,
but not individually rational w.r.t. the max-min value in pure strategies,
would not be a BCC equilibrium payoff,
so that a Nash equilibrium payoff need not be a BCC equilibrium payoff.

Moreover, in zero-sum games
it is not clear whether
there is a unique BCC equilibrium payoff.
If in zero-sum games there always is a unique BCC equilibrium payoff, then this quantity
can be called the {\em BCC value} of the game.
However, it is possible that in zero-sum games there will be more than one BCC equilibrium payoff,
in which case even in this class of games, the outcome will crucially depend on the relative
computational power the players have.

\subsubsection{A more general definition of a BCC equilibrium}

The definition of $c$-BCC equilibrium assumes that the utility of each player is additive,
and that the memory cost is linear in the memory size.
There are applications where the utility function $U_i$ has a different form.
\begin{itemize}
\item Players may disregard the memory cost, but be bounded by the size of memory that they use.
\[ U_i(M_1,M_2) = \left\{
\begin{array}{lll}
\gamma_i(M_1,M_2) & \ \ \ \ \ & |M_i| \leq k_i,\\
-\infty & & |M_i| > k_i.
\end{array}
\right. \]
This situation occurs, e.g., when players are willing to invest huge amount of money even if the profit is low,
but the available technology does not allow them to increase their memory size beyond some limit.
Such situation may occur, e.g., in the area of code breaking,
where countries invest large sums of money to be able to increase
the number of codes of other countries that they break,
and they are only bounded by technological advances.
\item
Memory is costly, yet players do not save money be reducing their memory size.
That is, a pair of mixed automata $(M_1,M_2)$ is a $c$-BCC equilibrium if for each $i \in \{1,2\}$
and for every pure automaton $P_i \in M_i$ one has
$\gamma_i(M_i,M_{3-i}) \geq \gamma_i(P_i,M_{3-i})$,
and, if $P_i > M_i$, one has
$\gamma_i(M_i,M_{3-i}) \geq \gamma_i(P_i,M_{3-i}) - c(|P_i| - |M_i|)$.
This situation occurs, e.g., when the players are organizations whose size cannot be reduced.
\end{itemize}

It may be of interest to study the set of equilibrium payoffs for various utility functions $U_i$,
and to see whether and how this set depends on the shape of this function.

\subsubsection{More than two players}

The concept of BCC equilibrium payoff is valid to games with any number of players.
However, Theorem \ref{th1} holds only for two-player games.
One crucial point in our construction is that if a deviation is detected,
a player is punished for a long (yet finite) period of time by a punishing action.
When there are more than two players,
the punishing action of, say, player 1 against player 2 may be different that
the punishing action of player 1 against player 3.
It is not clear how to construct an automaton that can punish each of the other players, if necessary,
and such that all these memory cells will be used on the equilibrium path.

\subsubsection{BCC equilibria in one-shot games}

The concept of BCC equilibrium that we presented here applies to repeated games.
However, the concept can be naturally adapted to one-shot games as well%
\footnote{We thank Ehud Kalai for drawing our attention to this issue.}.
For example, consider the following game, that appears in Halpern and Pass \cite{HP}.
Player 1 chooses an integer $n$ and tells it to player 2;
player 2 has to decide whether $n$ is a prime number or not, winning 1 if she is correct,
losing 1 if she is incorrect.
Plainly the value of this game is 1:
player 2 can check whether the choice of player 1 is a prime number.
However, as there is no efficient algorithm to check whether an integer is a prime number,
it is not clear whether in practice risk-neutral people would be willing to participate
in this game as player 2.

The concept of BCC equilibrium can be applied in such situations,
and one can study the set of BCC equilibrium payoffs, and how this set depends on the relative
memory cost of the two players.

In the context of the Computer Science literature one could conceive of an analog solution concept,
where automata are replaced by Turing machines,
and the memory size is replaced by the length of the machine's tape.

\section{BCC Equilibria in the Prisoner's Dilemma}
\label{section example}

In the present section we prove Theorem \ref{th1} for the Prisoner's Dilemma.
The construction in this case contains all the ingredients of the general case,
yet the simplicity of the Prisoner's Dilemma allows one to concentrate on the main aspects of the construction.
In Section \ref{section general} we indicate how to generalize this basic construction to general two-player repeated games.

The Prisoner's Dilemma is the two-player game depicted in Figure \arabic{figurecounter},
where each player has two actions: $A_1 = A_2 = \{Cooperate,Defect\}$.

\centerline{
\begin{picture}(140,85)(-60,0)
\put(-60, 18){Player 1}
\put(23, 68){Player 2}
\put(-10, 8){$C$}
\put(-10,28){$D$}
\put( 20,50){$D$}
\put( 60,50){$C$}
\put( 0,0){\numbercellong{$0,4$}{}}
\put( 0,20){\numbercellong{$1,1$}{}}
\put( 40,0){\numbercellong{$3,3$}{}}
\put( 40,20){\numbercellong{$4,0$}{}}
\end{picture}
}

\centerline{Figure \arabic{figurecounter}: The Prisoner's Dilemma.}
\addtocounter{figurecounter}{1}

\bigskip

The min-max level of each player is 1, and the punishing action of each player is $D$.
The set of feasible and (weakly) individually rational payoffs appear in Figure \arabic{figurecounter}.
It is equal to the quadrilateral $W$ with extreme points $(1,1)$, $(1,3\frac{2}{3})$, $(3,3)$ and $(3\frac{2}{3},1)$.


\centerline{Figure \arabic{figurecounter}: The feasible and individually rational payoffs in the Prisoner's Dilemma.}
\addtocounter{figurecounter}{1}

\bigskip

We now show that every feasible and individually rational payoff
vector $x$ is a BCC equilibrium payoff.
In the construction we do not use the special structure of the payoff matrix;
all we use is that each player has two actions, and that $D$ is the punishing action of both players.

Observe that each point in $W$ can be written as a convex combination of three vectors
in the payoff matrix, $(3,3)$, $(1,1)$, and either $(0,4)$ or $(4,0)$.
Assume w.l.o.g. that the latter holds, so that
\begin{equation}
x = \alpha_1(1,1) + \alpha_2(4,0) + \alpha_3(3,3),
\end{equation}
where $\alpha_1+\alpha_2+\alpha_3 = 1$ and $\alpha_1,\alpha_2,\alpha_3 \geq 0$.

Our goal is to define two sequences of mixed automata $(M_1(k))_k$ and $(M_2(k))_k$,
that support $x$ as a BCC equilibrium payoff:
the long-run average payoff under $(M_1(k),M_2(k))$ will converge to $x$.
The road-map of the proof is as follows.
We fix $k \in \dN$, and we define a play path $\omega^*$
that depends on $k$ and that will be the equilibrium path under $(M_1(k),M_2(k))$
(Section \ref{subsection equilibrium play}).
We then calculate a lower bound to the complexity of the play path for each player
(the complexity is of the order $k^3$, see Section \ref{section complexity}).
Recall that the complexity of a play path w.r.t. a player is the size of the smallest automaton
for that player that can implement this play path,
provided the other player follows her part in the play path.
We then construct, for each player,
a family of pure automata with this smallest size that implement the play path (Sections \ref{section1} and \ref{section2}).
We let the mixed automaton of each player choose one of these pure automata,
and finally we prove that each of these mixed automata is a $z(k)$-BCC best reply against the other,
where $\lim_{k \to \infty} z(k) = 0$ (see Sections \ref{section proof1} and \ref{section proof2}).

\subsection{The Equilibrium Play}
\label{subsection equilibrium play}

We fix throughout a natural number $k$, sufficiently large to satisfy several
conditions that will be set in the sequel.
Let $k_0$ be the largest integer that satisfies $(k_0)^2 + k_0 \leq k$.
We here define a specific play path $\omega^*$ that will be the equilibrium path.

We approximate $(\alpha_1,\alpha_2,\alpha_3)$ by rational numbers with denominator $k_0$;
that is, let $(k_1,k_2,k_3)$ be three natural numbers that satisfy
(a) $k_1+k_2+k_3=k_0$, and
(b) $\| \frac{k_j}{k_0} - \alpha_j\| \leq \frac{1}{\sqrt{k_0}}$ for $j=1,2,3$.
Let $k$ be a sufficiently large integer such that there are at least $k_2$ prime numbers larger than $k_2$ and smaller than ${k-k_1}$.
Because the number of prime numbers smaller than $k$ is approximately $\frac{k}{\ln(k)}$,
$k$ is of the order%
\footnote{In (b) we require that $\| \frac{k_j}{k_0} - \alpha_j\| \leq \frac{1}{\sqrt{k_0}}$
rather than $\| \frac{k_j}{k_0} - \alpha_j\| \leq \frac{1}{{k_0}}$, to accommodate the case $\alpha_3=0$.
If $\alpha_3=0$, with the latter requirement we would have $k_3 \in \{0,1\}$,
and there would not be $k_2$ prime numbers between $k_2$ and $k-k_1$.}
of $k_2\ln(k_2)$.

Let $\omega_0$ be the following play of length $k_0$ that generates a payoff close to $x$:
\begin{eqnarray}
\omega_0 &=& k_1 \times (D,D) + k_2 \times (D,C) + k_3 \times (C,C)\\
&=&
\begin{array}{ccc}
\underbrace{(D,D), \cdots, (D,D)}, & \underbrace{(D,C),\cdots,(D,C)}, &\underbrace{(C,C),\cdots,(C,C)}.\\
k_1 \hbox{ times} & k_2 \hbox{ times} & k_3 \hbox{ times}
\end{array}
\end{eqnarray}
Here, the notation $n \times a$ means a repetition of $n$ times the action pair $a$,
and $\omega_1 + \omega_2$ means the concatenation of $\omega_1$ and $\omega_2$.
Because of the choice of $(k_1,k_2,k_3)$, the average payoff along $\omega_0$ is $\frac{12}{\sqrt{k_0}}$-close to $x$.

Let $\omega^*$ be the play path that consists of the followings three parts:
\begin{itemize}
\item   A {\em punishment phase} that consists of $k^3$ times playing $(D,D)$.
\item   A {\em babbling phase}, that consists of $2k+1$ blocks:
in odd blocks (except the last one) the players play $k$ times $(C,C)$,
in even block they play $k$ times $(D,D)$,
and in the last block the players play $k+1$ times $(C,C)$.
\item   A {\em regular play}, in which the players repeatedly play $\omega_0$.
\end{itemize}
Formally, the play path $\omega^*$ is:
\[ \omega^* = \underbrace{k^3 \times (D,D)}_{\hbox{Punishment}} +
\underbrace{\sum_{n=1}^k \big(k\times(C,C) + k\times(D,D)\big) +
(k+1)\times(C,C)}_{\hbox{Babbling}} + \underbrace{\sum_{n=1}^\infty \omega_0}_{\hbox{Regular}}. \]

The roles of the three phases are as follows.
\begin{itemize}
\item   As in Abreu and Rubinstein  \cite{AR88}, the punishment phase ensures that punishment is on the equilibrium path.
Because the players minimize their automaton size,
subject to maximizing their payoff,
if the punishment phase was off the equilibrium path,
players could save states by not implementing it.
But if a player does not implement punishment, the other player may safely deviate, knowing that she will not be punished.
In our construction, detectable deviations of the other player
will lead the automaton to restart and re-implement $\omega^*$,
thereby initiating a long punishment phase.
The length of the punishment phase, $k^3$,
is much longer than the babbling phase to ensure that the punishment is severe.
\item  The importance of the babbling phase is that it allows us to build up the mixed strategy equilibrium.
To reach any equilibrium payoff in the convex hull, players need to implement sequences of action pairs. Some of them could be
played by means of  some previously used states. Nevertheless this construction may fail due to the possible deviation (without punishment) of the opponent.
In order to avoid this weakness, it must be concealed the position of such re-used states. It is here where the use of the mixed strategy plays a decisive role: to hide the chosen pure strategy. This set of pure strategies will be characterised by the location of the re-used states within a convenient set of states. In our construction this is implemented by the babbling phase.

The babbling phase which serves two purposes.
First, because it is coordinated, it is not difficult to calculate its complexity.
Second, it is sufficiently long, so that to implement the regular phase one does not need new states,
but rather one can re-use states that implement the babbling phase.
Moreover, its long lengths ensures that, if the states that are re-used are chosen randomly,
to find which states are re-used with non-negligible probability the other player must use a very large automaton:
to profit by deviating the other player needs to search for the re-used states,
a task that requires a significantly larger automaton than the one she currently uses.
\item   On the equilibrium path the regular play will be played repeatedly,
so that the long-run average payoff will be the average payoff along $\omega_0$, which is close to $x$.
\end{itemize}

\subsection{The complexity of $\omega^*$}
\label{section complexity}

Let $\omega$ be a (finite or infinite) sequence of action pairs.
We say that a mixed automaton $M_i$ of player $i$ is {\em compatible} with the play $\omega$ if,
when the other player $3-i$ plays her part in $\omega$,
the automaton generates the play of player $i$ in $\omega$ (with probability 1).
Plainly, different automata may be compatible with the same sequence $\omega$.
The {\em complexity} of $\omega$ w.r.t. player $i$ is the size of the smallest automaton of player $i$
that is compatible with $\omega$.
This concept was first defined and studied by Neyman \cite{N98},
who also provided a simple way to calculate it.

Our goal now is to calculate the complexity of $\omega^*$ w.r.t. the two players.

\begin{lemma}
\label{lemma 2}
The complexity of $\omega^*$ w.r.t. player 1 is $k^3 + 2k^2 + 1$,
and its complexity w.r.t. player 2 is $k^3 + 2k^2 + k+1$.
\end{lemma}

In the rest of this subsection we prove that the complexity of $\omega^*$ w.r.t.
each of the players is at least the quantities given in Lemma \ref{lemma 2}.
In the next two subsections we provide an automaton for player 1 (resp. for player 2) with size $k^3 + 2k^2 + 1$
(resp. $k^3 + 2k^2 + k+1$) that is compatible with $\omega^*$,
thereby completing the proof of Lemma \ref{lemma 2}.

We start by recalling Neyman's \cite{N98} characterization for the complexity of a play w.r.t. a player.

Denote by $\omega_t$ the sequence $\omega$ after deleting the first $t-1$ elements from the sequence.%
\footnote{If $\omega$ is a finite play, and $t$ is larger than the length of $\omega$,
then $\omega_t$ is an empty sequence of action pairs.}
Given a sequence of action pairs $\omega$, finite or infinite,
define an equivalence relation on the set of natural numbers $\dN$ as follows:
$t$ is equivalent (for player $i$) to $t'$ if any automaton of player $i$ that is compatible with $\omega_t$
is also compatible%
\footnote{In particular, the empty play is equivalent to any other play.}
with $\omega_{t'}$.
Denote this equivalence relation by $\sim_{\omega,i}$.
Neyman \cite{N98} proved that the complexity of $\omega$ w.r.t. a player is the number of equivalence classes
in this equivalence relation.

\subsubsection{The complexity of $\omega^*$ w.r.t. player 1 is at least $k^3 + 2k^2 + 1$}

The complexity of a sequence is at least the complexity of any of its subsequences (Lemma 2 in Neyman \cite{N98}).
To bound the complexity of $\omega^*$ w.r.t. player 1 we calculate the complexity w.r.t. player 1 of
the following prefix $\omega^*(1)$ of $\omega^*$:
\[ \omega^*(1) = k^3 \times (D,D) + \sum_{n=1}^k \big(k\times(C,C) + k\times(D,D)\big) + (k+1)\times(C,C). \]
In $\omega^*(1)$ the players play a coordinated play, i.e., there exists a one-to-one relationship between
the actions played by player 1 and the actions played by player 2:
 in every stage either both players play $C$ or both players play $D$.
Therefore, for every $t$,
any automaton of player 1 that is compatible with $\omega^*_t(1)$ can ignore the actions of player 2.
Consequently, an automaton of player 1 generates a deterministic sequence of actions. This implies that  if $t_1 < t_2$,
and $t_1$ and $t_2$ are equivalent (w.r.t. $\sim_{\omega^*(1),1}$),
then $\omega^*_{t_2}(1)$ is a prefix of $\omega^*_{t_1}(1)$.

Because a sequence of $k+1$ times $C$ appears only at the end of the sequence $\omega^*(1)$,
it follows that $\omega^*_{t_2}(1)$ is not a prefix of $\omega^*_{t_1}(1)$ whenever $t_1<t_2 \leq k^3 + 2k^2+1$.
In particular, the complexity of $\omega^*$ to player 1 is at least $k^3 + 2k^2+1$.

\subsubsection{The complexity of $\omega^*$ w.r.t. player 2 is at least $k^3 + 2k^2 + k+1$}

To bound the complexity of $\omega^*$ w.r.t. player 2, we
calculate the complexity for player 2 of the following prefix $\omega^*(2)$ of $\omega^*$:
\[ \omega^*(2) = k^3 \times (D,D) + \sum_{n=1}^k \big(k\times(C,C) + k\times(D,D)\big) +
(k+1)\times(C,C) + k_1 \times (D,D) + 1 \times (D,C). \]
Apart of the last action pair, the play path $\omega^*(2)$ consists of a coordinated play.
Hence, analogously to the analysis for player 1,
for every $t$, any automaton of player 2 that is compatible with $\omega^*_t(2)$ can ignore the actions of player 1.
We now count the number of equivalence classes of the relation $\sim_{\omega^*(2),2}$.
The sequence $1 \times (C,C) + k_1\times (D,D)+1 \times (C,C)$ appears along $\omega^*(2)$
only after $k^3+2k^2+k+1$ stages in $\omega^*(2)$.
It follows that the number of equivalence classes of $\sim_{\omega^*(2),2}$ is at least $k^3+2k^2+k+1$.
In particular, the complexity of $\omega^*$ to player 2 is at least $k^3+2k^2+k+1$.

\subsection{An automaton $M_1$ for player 1}
\label{section1}

In this section we define a family of pure automata for player 1,
all have size $k^3 + 2k^2 + 1$.
Each automaton in  the family is compatible with $\omega^{*}$.
This will prove that the complexity of $\omega^{*}$ w.r.t. player 1 is $k^3+2k^2+1$, as stated in Lemma 1.
In section \ref{Mixed1} we define a mixed automaton for player 1
that is supported by pure automata in this family and that will be part of the
$d$-BCC equilibrium for a proper $d > 0$.

The automata in the family are parameterized by two parameters:
an integer $j \in \{1,2,\ldots, k-1 \}$ and a set $H=\{h_1,h_2, \ldots, h_{k_2}\}$ of $k_2$ integers.
The range of $h_1,h_2, \ldots, h_{k_2}$ will be defined in step 3 below where they are used.

Given a pair $(j,H)$ we proceed to construct a pure automaton  $P_1^{j,H}$ for player 1.
For clarity of the exposition, the construction is divided into three steps.
We start in step 1 by defining transitions that implement the prefix of length $k^3+2k^2+1$ of $\omega^{*}$.
In step 2 we add transitions that implement the next $k+k_1$ action pairs in $\omega^{*}$,
and in step 3 we add transitions that implement the rest of $\omega^{*}$.
In step 1 we will use all the states of $P_1^{j,H}$.
In step 2 and 3 we will re-use states for implementing the rest of $\omega^{*}$.
The mixed automaton that we will define later will choose $j$ and $H$ randomly, to conceal the states that are re-used.

The size of the automaton $P_1^{j,H}$ that we construct is $k^3+2k^2+1$.
Denote its states by the integers $Q=\{1,2,\ldots,k^3+2k^2+1\}$, where $1$ is the initial state.

\subsubsection{Step 1: Implementing the prefix of $\omega^{*}$ of length $k^3+2k^2+1$.}

The prefix of length $k^3+2k^2+1$ of $\omega^{*}$ is:
$$
{\omega_1}=k^3 \times (D,D)+ \sum_{n=1}^{k} \big(k \times (C,C) + k \times (D,D)\big) +  (C,C).
$$

This play consists of the punishment phase followed by $k$ pairs of blocks,
each block is made of a $C$-block and a $D$-block (both of length $k$).
The length of $\omega_1$ is equal to the size of the automaton,
and therefore a naive implementation is to have one state for each action of player 1 in $\omega_1$:
state $q \in Q$ will implement the $q$'th action pair in $\omega_1$.
Formally, we divide $Q$ to three sets:
\begin{enumerate}
\item   $Q^P = \{1,2,\ldots,k^3\}$: this is the set of all states that implement the punishment phase.
\item   $Q^C=\bigcup_{n=0}^{k-1} \{ k^3+2nk+1, \ldots, k^3+2nk+k \} \cup \{k^3+2k^2+1\}$:
this is the set of states in all $C$-blocks.
\item   $Q^D= \bigcup_{n=0}^{k-1} \{k^3+2nk+k+1, \ldots, k^3+2nk+2k\} $:
this is the set of states in all $D$-blocks.
\end{enumerate}

The output function is:
\[ f(q) = \left\{
\begin{array}{lll}
D & \ \ \ \  & q \in Q^P \cup Q^D,\\
C & & q \in Q^C,
\end{array}
\right.
\]
and the transition function is
\[ g(q,f(q)) = q+1, \ \ \ 1 \leq q < k^3 + 2k^2+1. \]
Because the play in $\omega_1$ is coordinated,
the transition is defined only if player 2 complies with the desired play $\omega_1$.
Figure \arabic{figurecounter} illustrates the first step in the construction of the automaton $P_1^{j,H}$.
In this figure, the initial state is the dotted circle to the left,
the white squares correspond to states where the action is $D$,
and the black circles correspond to states where the action is $C$.

\begin{picture}(430,50)
\hspace{-4.5mm}\put (0,30)
{$\underbrace{\odot\hspace{-1.7mm}\rightarrow\hspace{-1.5mm}
{\square} \hspace{-1.7mm}\rightarrow\hspace{-1.5mm}
{\square}\hspace{-1.7mm}\rightarrow\hspace{-1.5mm}
{\square}\hspace{-1.7mm}\rightarrow\hspace{-1.5mm}
\ldots \hspace{-1.7mm}\rightarrow\hspace{-1.5mm}
{\square}}\hspace{-1.7mm}\rightarrow\hspace{-1.5mm}
\underbrace{{\bullet}\hspace{-1.7mm}\rightarrow\hspace{-1.5mm}
{\bullet}\hspace{-1.7mm}\rightarrow\hspace{-1.5mm}
\ldots \hspace{-1.7mm}\rightarrow\hspace{-1.5mm}
{\bullet}}\hspace{-1.7mm}\rightarrow\hspace{-1.5mm}
\underbrace{{\square}\hspace{-1.7mm}\rightarrow\hspace{-1.5mm}
{\square}\hspace{-1.7mm}\rightarrow\hspace{-1.5mm}
\ldots \hspace{-1.7mm}\rightarrow\hspace{-1.5mm}
{\square}}\hspace{-1.7mm}\rightarrow\hspace{-1.5mm}
\ldots \hspace{-1.7mm}\rightarrow\hspace{-1.5mm}
\underbrace{{\bullet}\hspace{-1.7mm}\rightarrow\hspace{-1.5mm}
{\bullet}\hspace{-1.7mm}\rightarrow\hspace{-1.5mm}
\ldots \hspace{-1.7mm}\rightarrow\hspace{-1.5mm}
{\bullet}}\hspace{-1.7mm}\rightarrow\hspace{-1.5mm}
\underbrace{{\square}\hspace{-1.7mm}\rightarrow\hspace{-1.5mm}
{\square}\hspace{-1.7mm}\rightarrow\hspace{-1.5mm}
\ldots \hspace{-1.7mm}\rightarrow\hspace{-1.5mm}
{\square}}\hspace{-1.7mm}\rightarrow\hspace{-1.5mm}
{\bullet}
$}

\put (25,7) {Punishment}
\put (123,7) {$C$-block}
\put (200,7) {$D$-block}
\put (301,7) {$C$-block}
\put (375,7) {$D$-block}
\put (439,7) {$C$}
\end{picture}

\begin{center}
Figure \arabic{figurecounter}: An implementation of $\omega_1$.
\end{center}
\addtocounter{figurecounter}{1}

\subsubsection{Step 2: Implementing the next $k+k_1$ action pairs.}

We now add to the automaton $P_1^{j,H}$ transitions that implement the next $k+k_1$ action pairs in $\omega^{*}$, which are
$$
\omega_2=k \times (C,C)+ k_1 \times (D,D)
$$
Here we use the parameter $j$.
Because (a) the play $\omega_2$ starts with $k \times (C,C)$,
and (b) each $C$-block has length $k$ and is followed by a $D$-block
whose length is more than $k_1$, we can use the $j$'th $C$-block and the following $D$-block to implement $\omega_2$.
Therefore, to implement $\omega_2$ it is sufficient to add one transition to $P_1^{j,H}$,
from the last state to the beginning of the $j$'th $C$-block:
\[ g(k^3+2k^2+1,C) = k^3 + 2(j-1)k+1. \]
Figure \arabic{figurecounter} illustrate the automaton $P_1^{j,H}$ with this additional transition.

\vskip1cm
\begin{picture}(430,80)
\hspace{-4.5mm}\put (0,40)
{$\underbrace{\odot\hspace{-1.7mm}\rightarrow\hspace{-1.5mm}
{\square} \hspace{-1.7mm}\rightarrow\hspace{-1.5mm}
{\square}\hspace{-1.7mm}\rightarrow\hspace{-1.5mm}
{\square}\hspace{-1.7mm}\rightarrow\hspace{-1.5mm}
\ldots \hspace{-1.7mm}\rightarrow\hspace{-1.5mm}
{\square}}\hspace{-1.7mm}\rightarrow\hspace{-1.5mm}
\underbrace{{\bullet}\hspace{-1.7mm}\Rightarrow\hspace{-1.5mm}
{\bullet}\hspace{-1.7mm}\Rightarrow\hspace{-1.5mm}
\ldots \hspace{-1.7mm}\Rightarrow\hspace{-1.5mm}
{\bullet}}\hspace{-1.7mm}\Rightarrow\hspace{-1.5mm}
\underbrace{{\square}\hspace{-1.7mm}\rightarrow\hspace{-1.5mm}
{\square}\hspace{-1.7mm}\rightarrow\hspace{-1.5mm}
\ldots \hspace{-1.7mm}\rightarrow\hspace{-1.5mm}
{\square}}\hspace{-1.7mm}\rightarrow\hspace{-1.5mm}
\ldots \hspace{-1.7mm}\rightarrow\hspace{-1.5mm}
\underbrace{{\bullet}\hspace{-1.7mm}\rightarrow\hspace{-1.5mm}
{\bullet}\hspace{-1.7mm}\rightarrow\hspace{-1.5mm}
\ldots \hspace{-1.7mm}\rightarrow\hspace{-1.5mm}
{\bullet}}\hspace{-1.7mm}\rightarrow\hspace{-1.5mm}
\underbrace{{\square}\hspace{-1.7mm}\rightarrow\hspace{-1.5mm}
{\square}\hspace{-1.7mm}\rightarrow\hspace{-1.5mm}
\ldots \hspace{-1.7mm}\rightarrow\hspace{-1.5mm}
{\square}}\hspace{-1.7mm}\rightarrow\hspace{-1.5mm}
{\bullet}
$}
\put (444,47){\line(0,1){30}}
\put (444,77){\line(-1,0){324}}
\put (120,77){\vector(0,-1){30}}
\put (25,12) {Punishment}
\put (123,12) {$C$-block}
\put (200,12) {$D$-block}
\put (301,12) {$C$-block}
\put (375,12) {$D$-block}
\put (439,12) {$C$}
\end{picture}

\begin{center}
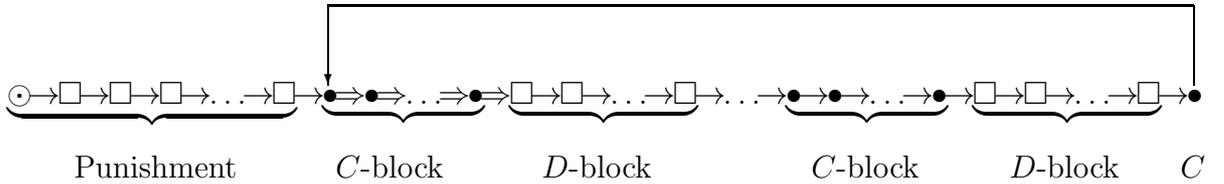
Figure \arabic{figurecounter}: The automaton $P_1^{j,H}$ after the second step.
\end{center}
\addtocounter{figurecounter}{1}

\subsubsection{Step 3: Implementing the rest of $\omega^*$.}

We now add to the automaton $P_1^{j,H}$ transitions that implement the next $k_2+k_3$ action pairs in $\omega^{*}$,
which are
 $$
 \omega_3=k_2 \times (D,C) + k_3 \times (C,C),
 $$
and continue to implement the regular play, which is a periodic repetition of $\omega_0$.

Here we use the parameter set $H$.
To implement the $k_2$ repetitions of $(D,C)$ we re-use states in a $D$-block,
whose identity is determined by the set $H$.
Thus, whenever in a re-used state,
if player 2 plays $D$, the automaton $P_1^{j,H}$ assumes that the play is in the babbling phase,
whereas if player 2 plays $C$, the automaton assumes that $\omega_3$ is implemented.
Because $\omega_3$ comes after a sequence $k_1 \times (D,D)$, the first re-used state must be the
$k_1+1$ state in the $j$'th $D$-block.
Because after the sequence $k_3 \times (C,C)$ the play continues with the next repetition of $\omega_0$,
namely, with $k_1 \times (D,D)$,
the sequence $k_3 \times (C,C)$ will be implemented at the end of the $j$'th $C$-block.

Formally, assume that the set $H$ satisfies the following two conditions:
\label{condition D}
\begin{enumerate}
\item[(D1)] $h_1=k^3+2(j-1)k+k+k_1+1$, and
\item[(D2)] $h_2,h_3, \ldots,h_{k_2}$ are distinct states in $Q^D$, all different from $h_1$.
\end{enumerate}
We add the following transitions (see Figure \arabic{figurecounter}):
\begin{eqnarray}
g(h_n,C) &=& h_{n+1}, \ \ \ 1 \leq n < k_2-1,\\
g(h_{k_2},C) &=& k^3 + 2(j-1)k + (k-k_3).
\end{eqnarray}
In Figure \arabic{figurecounter}, re-used states are denoted by triangles.
When the automaton $P_1^{j,H}$ is at such a state
it plays the action $D$;
if player 2 plays the action $D$, the transition is to the subsequent (square) state,
whereas if player 2 plays $C$, the transition is to the next triangle state.

\vskip1cm
\begin{picture}(430,80)
\hspace{-4.5mm}\put (0,40)
{$\hspace{-2.9mm}
\ldots{\bullet}\hspace{-1.7mm}\Rightarrow\hspace{-1.5mm}
{\bullet}\hspace{-1.7mm}\Rightarrow\hspace{-1.5mm}
\ldots \hspace{-1.7mm}\Rightarrow\hspace{-1.5mm}
{\bullet}\hspace{-1.7mm}\Rightarrow\hspace{-1.5mm}
\underbrace{{\bullet}\hspace{-1.7mm}\Rightarrow\hspace{-1.5mm}
{\bullet}\hspace{-1.7mm}\Rightarrow\hspace{-1.5mm}
{\bullet}}\hspace{-1.7mm}\Rightarrow\hspace{-1.5mm}
\underbrace{{\square}\hspace{-1.7mm}\Rightarrow\hspace{-1.5mm}
{\square}\hspace{-1.7mm}\Rightarrow\hspace{-1.5mm}
{\square}}\hspace{-1.7mm}\Rightarrow\hspace{-1.5mm}
\underbrace{{\triangle}\hspace{-1.7mm}\rightarrow\hspace{-1.5mm}
{\square}\hspace{-1.7mm}\rightarrow\hspace{-1.5mm}
{\triangle}\hspace{-1.7mm}\rightarrow\hspace{-1.5mm}
{\square}\hspace{-1.7mm}\rightarrow\hspace{-1.5mm}
{\triangle}\hspace{-1.7mm}\rightarrow\hspace{-1.5mm}
{\square}\hspace{-1.7mm}\rightarrow\hspace{-1.5mm}
{\triangle}}\hspace{-1.7mm}\rightarrow\hspace{-1.5mm}
{\square}\hspace{-1.7mm}\rightarrow\hspace{-1.5mm}
\ldots \hspace{-1.7mm}\rightarrow\hspace{-1.5mm}
{\square}\hspace{-1.7mm}\rightarrow\hspace{-1.5mm}
\ldots \hspace{-1.7mm}\rightarrow\hspace{-1.5mm}
{\square}\hspace{-1.7mm}\rightarrow\hspace{-1.5mm}
\ldots \hspace{-1.7mm}\rightarrow\hspace{-1.5mm}
{\square}\hspace{-1.7mm}\rightarrow\hspace{-1.5mm}
{\bullet}
$}
\put (190,45){\line(0,1){15}}
\put (190,60){\line(1,0){25}}
\put (215,60){\vector(0,-1){15}}
\put (228,45){\line(0,1){15}}
\put (228,60){\line(1,0){25}}
\put (253,60){\vector(0,-1){15}}
\put (266,45){\line(0,1){15}}
\put (266,60){\line(1,0){25}}
\put (291,60){\vector(0,-1){15}}
\put (305,40){\line(0,-1){15}}
\put (305,25){\line(-1,0){228}}
\put (77,25){\vector(0,1){15}}
\put (467,47){\line(0,1){30}}
\put (467,77){\line(-1,0){457}}
\put (10,77){\vector(0,-1){30}}
{\small
\put (70,0) {$k_3 \times (C,C)$}
\put (130,0) {$k_1 \times (D,D)$}
\put (215,0) {$k_2 \times (D,C)$}
}
\end{picture}

\begin{center}
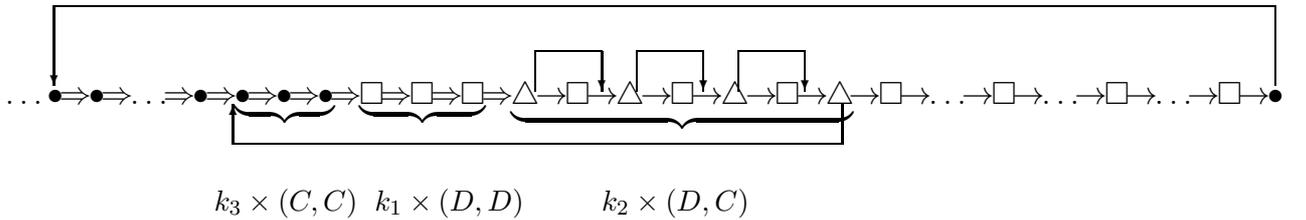
Figure \arabic{figurecounter}: The $j$'th $C$-block and $D$-block in $P_1^{j,H}$.
\end{center}
\addtocounter{figurecounter}{1}

\subsubsection{Last step: Deviations.}

By construction, the automaton $P_1^{j,H}$ is compatible with $\omega^*$.
In particular, the complexity of $\omega^*$ w.r.t. player 1 is $k^3+2k^2+1$ as stated in Lemma \ref{lemma 2}.
We now add transitions to detect deviations of player 2 as follows:
all transitions that were not defined in steps 1-3 lead to state 1.

Only re-used states accept both actions of player 2;
the other states accept only the action that is indicated by $\omega^*$.
Because a punishment phase of length $k^3$ begins in state 1,
any deviation in a non re-used state is followed by a severe punishment.
In the next subsection we define the mixed automaton that player 1 uses.
The parameters $j$ and $H$ will be chosen randomly,
so that to profit by deviation, player 2 will have to learn $j$ or $H$,
and such a learning process requires a large memory.

\subsubsection{Mixed strategy}\label{Mixed1}

We now define the mixed automaton $M_1=M_1(k)$ for player 1.
For every $n$, $1 \leq n \leq k_2$, define
\[ \widehat H_n = \{1, 1+n, 1+2n, 1+3n, \ldots, 1+k_2n\}, \]
and
\[ H_n = \{ k^3 + 2(n-1)k + k + k_1 + h \colon h \in \widehat H_n\}. \]
Thus, $H_n$ contains $k_2$ states in the $n$'th $D$-block, that are equally spaced,
and the distance between each two adjacent states is $n$.
Because $k_0 + (k_0)^2 \leq k$,
there are enough states in the $D$-block to accommodate this construction,
and the two conditions (D1) and (D2) (in page \pageref{condition D}) are satisfied.

Let $J = \{j_1,j_2,\ldots,j_{k_2}\}$ be a collection of $k_2$ different prime numbers in the range $\{k_2+1,k_2+2,\ldots,{k-k_1}\}$,
which exist by the choice of $k$.
Let $M_1$ be the mixed automaton of player 1 that assigns a probability $\frac{1}{k_2}$
to each of the pure automata $P_1^{j_l,H_{j_l}}$.

\subsection{An automaton $M_2$ for player 2} \label{section2}

In this section we describe an analog construction to the one we presented in section \ref{section1},
for a mixed automaton of player 2.
We construct a family of pure automata for player 2,
all of size $k^3 + 2k^2+k+1$, and all compatible with $\omega^*$ for player 2.
As for player 1, the automata in the family depend on two parameters,
an integer $j \in \{1,2,\ldots,k-1\}$ and a set $H$ of integers.

\subsubsection{Step 1: Implementing the prefix of length $k^3+2k^2+k+1$ of $\omega^*$.}

We start by implementing the prefix
of length $k^3+2k^2+k+1$ of $\omega^{*}$ by a naive automaton with $k^3+2k^2+k+1$ states.
The prefix is:
$$
\omega_4=k^3 \times (D,D)+ \sum_{n=0}^{k-1} \Large(k \times (C,C) + k \times (D,D)\Large) +  (k+1) \times (C,C),
$$
and it contains the punishment phase and the babbling phase.
As for player 1, we define an automaton that implements each action pair in one state.
\smallskip
Let $Q = \{1,2,\ldots,k^3+2k^2+k+1\}$ be the set of states of the automaton, and
divide $Q$ into three sets, as follows:
\begin{enumerate}
\item   $Q^P = \{1,2,\ldots,k^3\}$: this is the set of all states that implement the punishment phase.
\item   $Q^C=\bigcup_{n=0}^{k-1} \{ k^3+2nk+1, \ldots, k^3+2nk+k \} \cup \{k^3+2k^2+1,\cdots, k^3+2k^2+k+1\}$:
this is the set of all states in $C$-blocks.
\item   $Q^D= \bigcup_{n=0}^{k-1} \{k^3+2nk+k+1, \ldots, k^3+2nk+2k\} $:
this is the set of all states in $D$-blocks.
\end{enumerate}

The output function is:
\[ f(q) = \left\{
\begin{array}{lll}
D & \ \ \ \  & q \in Q^P \cup Q^D,\\
C & & q \in Q^C,
\end{array}
\right.
\]
and the transition function is
\[ g(q,f(q)) = q+1, \ \ \ 1 \leq q < k^3 + 2k^2+k+1. \]

\subsubsection {Step 2: Implementing the next $k_1$ actions in $\omega^*$.}

We now add the transitions that implement the next $k_1$ actions in $\omega^*$,
which are $\omega_5 = k_1 \times (D,D)$.
To this end, we re-use states in a $D$-block, and because after $\omega_5$ player 2 plays $C$ in $\omega^*$,
we re-use the last $k_1$ states that implement a $D$-block.
So that player 1 does not know which $D$-block is re-used, we use the $j$'th $D$-block.
Formally,
\[ g(k^3+2k^2+k+1,C)=k^3+2k j-k_1. \]

\subsubsection {Step 3: Implementing the rest of $\omega^*$.}

We now add transitions that implement $\omega_6=k_2 \times(D,C) +k_3 \times(C,C)+\sum_{n=1}^\infty \omega_0$.
To implement the sequence $k_2 \times(D,C)$ we
re-use states in a $C$-block that are determined by the set $H = \{h_1,h_2,\ldots,h_{k_2}\}$.
The first re-used state must be the first state in the $j+1$'th $C$-block,
and therefore
$h_1 = k^3+2k j+1$.
Because the second part of $\omega_3$, that is, $k_3 \times(C,C)$, should lead to the sequence $k_1 \times (D,D)$
that starts $\omega_0$,
the states that implement that part must be the last $k_3$ states in $Q$;
therefore we must have
$h_{k_2} = k^3+2k^2+k+1-k_3$.
Finally, we require that
$h_1,h_2,h_3,\ldots,h_{k_2}$ are distinct states in $C$-blocks.

Transitions are defined as follows:
\begin{eqnarray}
g(h_n, D) &=& h_{n+1}, \ \ \ \ \ 1 \leq n < k_2,\\
g(h_{k_2},D) &=& h_{k_2}+1.
\end{eqnarray}

\subsubsection{Last step: Deviations.}

Finally we add transitions to handle deviations in states that are not re-used.
All transitions that are not defined in steps 1-3, lead to state 1,
so that such deviations initiate a long punishment phase.

\subsubsection{Mixed strategy of player 2.}

The definition of the mixed strategy $M_2=M_2(k)$ is analog to that of $M_1$.
Recall that
\[ \widehat H_n = \{1, 1+n, 1+2n, 1+3n, \ldots, 1+k_2n\}, \]
and
\[ H_n = \{ k^3 + 2nk + k + k_1 + h \colon h \in \widehat H_n\}, \]
and that $J$ is a set of $k_2$ distinct prime numbers in the range $\{k_2+1,\ldots,{k-k_1}\}$.

Let $M_2$ be the mixed automaton of player 1 that assigns a probability $\frac{1}{k_2}$
to each of the pure automata $P_2^{j_l,H_{j_l}}$.

In the following subsections we show that the sequence $(M_1(k),M_2(k))_k$ supports $x$ as a BCC equilibrium payoff.
That is, we show that (1) the expected long-run average payoff under $(M_1,M_2)$ is $\frac{12}{\sqrt{k_0}}$-close to $x$,
(2) no player can profit by deviating to a smaller automaton,
and (3) we bound the amount a player can profit by deviation to a larger automaton.

\subsection{The expected payoff under $(M_1,M_2)$ is close to $x$}
\label{section proof1}

By construction, the automaton $M_1$ (resp. $M_2$) is compatible with the play $\omega^*$ for player 1 (resp. player 2).
Therefore,
if the players use these automata the play is $\omega^*$,
and the long-run average payoff is $\frac{12}{\sqrt{k_0}}$-close to $x$.

Define
\[ x^* := \gamma(M_1,M_2) =
\frac{k_1}{k_1+k_2+k_3}u(D,D) + \frac{k_2}{k_1+k_2+k_3}u(D,C) + \frac{k_3}{k_1+k_2+k_3}u(C,C). \]

\subsection{$(M_1,M_2)$ is a $c$-BCC Equilibrium}
\label{section proof2}

In this section we prove that
$(M_1,M_2)$ is a $c$-BCC-equilibrium, for every $c$ that satisfies $\frac{3}{k_2\times k^3} < c < \frac{\eta}{2k^3}$.
We only prove the claims for player 2.
The claims for player 1 can be proven analogously.
Below we denote the state of an automaton of player $i$ at stage $t$ by $q_i(t)$.

For $l \in \{1,2,\ldots,k_2\}$ denote $P_1^l := P_1^{j_l,H_l}$,
so that the support of $M_1$ is
$P_1^1,P_1^2,\ldots,P_1^{k_2}$.
Let $j^l$ and $H^l$ be the parameters $j$ and $H$ of $P_1^l$, for $l=1,2,\ldots,k_2$.
Let $P_2$ be an arbitrary pure automaton that implements a strategy of player 2.
We denote by $\omega^l$ the play that is generated under $(P_1^l,P_2)$.

Suppose that the players use the automata $(P_1^l,P_2)$.
If $P_2$ is not compatible with $\omega^*$ for player 2,
then $P_1^l$ restarts whenever a deviation from $\omega^*$ is detected,
and a punishment phase starts.
Denote by $t_n^l$ the stage at the $n$'th time in which $P_1^l$ visits state 1 when facing $P_2$:
\begin{eqnarray*}
t_1^l &:=& 1,\\
t_{n+1}^l &:=& \min\left\{t > t_n^l \colon q_1(t) = 1\right\}, \ \ \ \ \ n \geq 1.
\end{eqnarray*}
By convention, the minimum of an empty set is $\infty$.

There are two scenarios where player 2 may improve her long-run average payoff.
One possibility is if there exists $n$ such that $t_n^l < \infty = t_{n+1}^l$.
Then $t_n^l$ is the last stage in which the automaton $P_1^l$ restarts;
in other words, this is the last stage in which a punishment phase starts.
If the play after stage $t_n^l$ is different than $\omega^*$,
it means that player 2 plays as if she knows $j^l$ and/or $H^l$,
and she might use this information to improve her payoff.
Another possibility is
that $(t_n^l)_{n\in \dN}$ are finite and between two of
these stages the average payoff of player 2 is higher than $x^*_2$
(in fact, if $(t_n^l)_{n\in \dN}$ are finite then, so that player 2 improves her payoff,
the average payoff between $t_n^l$ and $t_{n+1}^l-1$ should be higher than $x^*$ infinitely often).

This leads us to the following definition.

\begin{definition}
\label{def fool}
The automaton $P_2$ {\em fools} the automaton $P_1^l$ if either one of the following conditions hold:
\begin{enumerate}
\item[C1)]   There is $n_0 \in \dN$ such that $t_{n_0}^l<\infty=t_{n_0+1}^l$ and $\omega^l_{n_0} \neq \omega^*$.
\item[C2)]   $t_n^l < \infty$ for every $n \in \dN$,
and there is $n_0 \in \dN$ such that the average payoff for player 2
between stages $t_{n_0}^l$ and $t_{{n_0}+1}^l-1$ is strictly higher%
\footnote{Observe that in this case $t_{{n_0}+1}^l \geq t_{n_0}^l + k^3$.
In fact, a stronger bound can be obtained.}
than $x^*_2$.
\end{enumerate}
\end{definition}

If condition C1 holds, we say that $P_2$ fools $P_1^l$ in stages $\{t_{n_0}^l, t_{n_0}^l+1,\ldots\}$.
If condition C2 holds, we say that $P_2$ fools $P_1^l$ in stages $\{t_{n_0}^l, t_{n_0}^l +1,\ldots,t_{n_0+1}^l-1\}$.
In both cases we set $t_*^l = t_{n_0}^l$,
and we say that at stage $t_*^l$ player 2 starts to fool $P_1^l$.
Denote by $R_l = \{q_2(t_*^l), q_2(t_*^l+1),\cdots,q_2(t_*^l+k^3-1)\}$
the $k^3$ states that $P_2$ visits at the beginning of the period in which it fools $P_1^l$.
We will prove below that the sets $(R_l)_{l=1}^n$ are disjoint,
thereby bounding from below the size of any automaton of player 2 that obtains high payoff when facing $M_1$.

Neither C1 nor C2 imply that the long-run average payoff under $(P_1^l,P_2)$ is higher than $x^*_2$.
Yet, as the next lemma shows, the converse is true:
if the long-run average payoff of player 2 under $(P_1^l,P_2)$ exceeds $x^*_2$,
then $P_2$ must have fooled $P_1^l$.

\begin{lemma}
\label{lemma1}
If $P_2$ does not fool $P_1^l$ then $\gamma_2(P_1^l,P_2) \leq x^*_2$.
\end{lemma}

\begin{proof}
Since both $P_1^l$ and $P_2$ are automata, the long-run average payoff of player 2 under $(P_1^l,P_2)$ exists.
Suppose first that $t_n^l < \infty$ for every $n \in \dN$.
Because $P_2$ does not fool $P_1^l$,
for every $n \in \dN$ the average payoff of player 2 between stages $t_n^l$ and $t_{n+1}^{l}$ is at most $x^*_2$,
and therefore $\gamma_2(P_1^l,P_2) \leq x^*_2$.

Suppose now that there is $n_0 \in \dN$ such that $t_{n_0}^l<\infty=t_{{n_0}+1}^l$.
Because $P_2$ does not fool $P_1^l$, we have $\omega^l_{n_0}=\omega^*$,
so that the long-run average payoff of player 2 after stage $t_{n_0}^l$ is $x^*_2$, and the result follows.
\end{proof}

\bigskip

Our goal is to relate the number of pure automata $P_1^l$ that $P_2$ fools to the size of $P_2$.
In fact, we will prove that the size of $P_2$ is at least $k^3$ times the number of
pure automata $P_1^l$ that $P_2$ fools.
To this end, we now prove that if $P_2$ fools both $P_1^{l_1}$ and $P_1^{l_2}$,
then $R_{l_1}$ and $R_{l_2}$ are disjoint:
the automaton of player 2 uses different states to fool each of the two automata.

\begin{lemma}
\label{lemma4}
Let $1 \leq l_1 < l_2 \leq k_2$.
If $P_2$ fools both $P_1^{l_1}$ and $P_1^{l_2}$,
then $R_{l_1} \cap R_{l_2} =\emptyset$.
\end{lemma}

The subtle definition of $(j_l,H_l)_{l=1}^{k_2}$ is the key ingredient in the proof of Lemma \ref{lemma4}.
An immediate corollary of Lemma \ref{lemma4} is:
\begin{corollary}
\label{corollary5}
Denote by $L_0$ the number of pure automata $P_1^l$ that $P_2$ fools.
Then $|P_2| \geq L_0k^3$.
\end{corollary}

\noindent
\begin{proof}[Proof of Lemma \ref{lemma4}]

\noindent\textbf{Step 1:} If $P_2$ fools $P_1^l$ then the states in $R_l$ are distinct: $|R_l| = k^3$.

At stage $t_*^l$ the automaton $P_1^l$ restarts; it expects the sequence $k^3 \times (D,D) + k \times (C,C)$,
and none of the states $\{1,2,\ldots,k^3+1\}$ of $P_1^l$ is re-used.
Because this play is coordinated,
its complexity is $k^3+1$,
and therefore player 2 must use at least $k^3$ distinct states to implement its prefix of length $k^3$.

\noindent\textbf{Step 2:} If $R_{l_1}$ and $R_{l_2}$ are not disjoint, then
$q_2(t_*^{l_1}+k^3-1) = q_2(t_*^{l_2}+k^3-1)$: the last state in $R_{l_1}$ coincides with the last state in $R_{l_2}$.

Suppose that $R_{l_1}$ and $R_{l_2}$ are not disjoint,
and assume that $q_2(t_*^{l_1}+n_1) = q_2(t_*^{l_2}+n_2)$.
We argue that necessarily $n_1=n_2$.
Indeed, assume to the contrary that $n_1 < n_2$.
Because in the $k^3$ stages that follow stage $t_*^{l_1}$ the automaton $P_1^{l_1}$ plays $D$,
and in the $k^3$ stages that follow stage $t_*^{l_2}$ the automaton $P_1^{l_2}$ plays $D$,
the automaton $P_2$ receives the same inputs (when facing $P_1^{l_1}$ after stage $t_*^{l_1}$,
and when facing $P_1^{l_2}$ after stage $t_*^{l_2}$),
so that it evolves in the same way:
$q_2(t_*^{l_1}+n_1+s) = q_2(t_*^{l_2}+n_2+s)$ for every $s$ that satisfies $1 \leq s \leq k^3-n_2$.
Because $P_2$ fools $P_1^{l_1}$, the action $P_2$ plays in state $q_2(t_*^{l_1}+n_1+k^3-n_2+1)$ is $D$.
Because $P_2$ fools $P_1^{l_2}$, the action $P_2$ plays in state $q_2(t_*^{l_2}+n_2+k^3-n_2+1)$ is $C$.
But $q_2(t_*^{l_1}+n_1+k^3-n_2+1)=q_2(t_*^{l_2}+n_2+k^3-n_2+1)$, a contradiction.

Because in the first $k^3$ stages after visiting stage 1,
both $P_1^{l_1}$ and $P_1^{l_2}$ play in the same manner
(both output $D$), it follows that the evolution of $P_2$ when facing either $P_1^{l_1}$ or $P_1^{l_2}$ is the same.
The claim follows.

\noindent\textbf{Step 3:} $R_{l_1} \cap R_{l_2} = \emptyset$.

Assume to the contrary that $R_{l_1}$ and $R_{l_2}$ are not disjoint.
Denote by $(j_1,H_1)$ and $(j_2,H_2)$ the parameters $(j,H)$ of $P_1^{l_1}$ and $P_1^{l_2}$ respectively.
By Step 2, the last state in $R_{l_1}$ coincides with the last state in $R_{l_2}$.
Both automata $P_1^{l_1}$ and $P_1^{l_2}$ continue in the same way, until one of them observes a deviation,
in which case it restarts.

Denote by $t_*^{l_1}+n$ the first stage after stage $t_*^{l_1}$ in which $P_2$ deviates from $\omega^*$
when facing $P_1^{l_1}$.
Because $R_{l_1} = R_{l_2}$, the first stage after stage $t_*^{l_2}$ in which $P_2$ deviates from $\omega^*$
when facing $P_1^{l_2}$ is $t_*^{l_2}+n$.
Because $P_2$ fools both $P_1^{l_1}$ and $P_1^{l_2}$,
the state that $P_1^{l_1}$ visits in stage $t_*^{l_1}+n$ is a re-used state,
as is the state that $P_1^{l_2}$ visits in stage $t_*^{l_2}+n$.
Because the re-used states in $P_1^{l_1}$ are in the $j_1$'th $D$-block,
while the re-used states in $P_1^{l_1}$ are in the $j_2$'th $D$-block,
the automata $P_1^{l_1}$ and $P_1^{l_2}$ are both in re-used states only when they implement the action pairs $(D,C)$,
that is, in the second part of the regular play $\omega_0$.

Let us now verify that $P_2$ cannot fool both $P_1^{l_1}$ and $P_1^{l_2}$.
Because in a $D$-block both automata $P_1^{l_1}$ and $P_1^{l_2}$ play $D$ unless a deviation is detected
and a punishment phase starts,
the evolution of $P_2$, when facing either $P_1^{l_1}$ or $P_1^{l_2}$ is the same,
as long as these automata are in the $D$-block.
It is therefore sufficient to show that there is no sequence of actions of player 2
that differ from the play of $\omega^*$ in $D$-blocks,
and that does not initiate a punishment phase when facing either $P_1^{l_1}$ or $P_1^{l_2}$.

Because $H_1$ (resp. $H_2$) contains $k_2$ numbers, equally spaced with distance $j_1$ (resp. $j_2$),
the difference $h_{r_1}-h_{r_2}$ of pairs of elements in $H_1$ (resp. $H_2$)
is a multiple of $j_1$ (resp. $j_2$).
Because $j_1$ and $j_2$ are prime numbers larger than $k_2$,the differences generated by $H_1$ are different than those generated by $j_2$.
It follows that the unique two sequences of actions of player 2 that does not initiate a punishment phase
neither when facing $P_1^{l_1}$ in block $j_1$ nor when facing $P_1^{l_2}$ in block $j_2$ are
(a) a repetition of $k_2$ times $C$, and
(b) a repetition of $k-k_1$ times $D$.
Because $P_2$ deviates from $\omega^*$, only the sequence in (b) should be considered.

Now, in all blocks after block $j_1$ (resp. $j_2$), the automaton $P_1^{l_1}$ (resp. $P_1^{l_2}$) does not re-use states.
Because $P_2$ fools both $P_1^{l_1}$ and $P_1^{l_2}$, it must follow the play indicated by these automata.
However, because $j_1 \neq j_2$, when the first of these two automata reaches its last state,
that automaton initiates a punishment phase if $P_2$ plays $D$,
while the other initiates a punishment phase if $P_2$ plays $C$.
This implies that if $P_2$ plays the sequence in (b),
then it cannot fool both $P_1^{l_1}$ and $P_1^{l_2}$, as desired.
\end{proof}

\bigskip

Recall that the min-max value in pure strategies of both players is 1.
Therefore, $\min\{x_1^*-1,x_2^*-1\} > 0$ is the minimal difference between the target payoff $x^*$ and
the min-max value.
We now prove that player 2 cannot profit by deviating to an automaton smaller than $M_2$.
\begin{lemma}
Let $\eta < x^*_2-1$, and assume that $k$ is sufficiently large so that $\frac{4}{k_2} + \frac{8}{k} < \frac{\eta}{2}$.
Let $P'_2$ be an automaton for player 2 with size smaller than $k^3+2k^2+k+1$.
Then $\gamma_2(M_1,P'_2) -c|P'_2|\leq \gamma_2(M_1,M_2)- c|M_2|$,
provided $c < \frac{\eta}{2k^3}$.
\end{lemma}

\begin{proof}
Because the complexity of $\omega^*$ w.r.t. player 2 is $k^3+2k^2+k+1$,
the play under $(P_1^l,P'_2)$ is not $\omega^*$.
By Lemma \ref{lemma4}, and because the size of $P_2$ is smaller than $2k^3$,
the automaton $P'_2$ can fool at most one of the automata $(P_1^l)$.
Because it cannot generate $\omega^*$,
any automaton which $P_2$ does not fool restarts after at most $k^3+2k^2+k$ stages,
so that the average payoff is at most $\frac{k^3}{k^3+2k^2+k} + 4\frac{2k^2+k}{k^3+2k^2+k}$.
It follows that the expected payoff $\gamma_2(M_1,P'_2)$ is at most
\[ 4\frac{1}{k_2} + 4\frac{k_2-1}{k_2}\frac{2k^2+k}{k^3+2k^2+k} + \frac{k_2-1}{k_2}\frac{k^3}{k^3+2k^2+k}
\leq 1 + \frac{4}{k_2} + \frac{8}{k}. \]
Because the size of the automaton $M_2$ is $k^3+2k^2+k+1$,
the gain of reducing the size of automaton from $|M_2|$ to $|P'_2|$ is at most $c(k^3+2k^2+k)$.
So that player 2 does profit by this deviation, we need to require that
\[ x_2^* \geq 1 + \frac{4}{k_2} + \frac{8}{k} + c(k^3+2k^2+k),\]
and therefore it is enough to require that
\[ x_2^* - 1 > \eta > \frac{4}{k_2} + \frac{8}{k} + c(k^3+2k^2+k). \]
The right-hand side inequality holds provided
\[ c < \frac{\eta - \frac{4}{k_2} - \frac{8}{k}}{k^3+2k^2+k}, \]
so it is enough to require that $c < \frac{\eta}{2k^3}$.
\end{proof}

We finally prove that player 2 cannot profit by deviating to an automaton larger than $M_2$.

\begin{lemma}
\label{lemma6}
Let $P'_2$ be a pure automaton such that $\gamma_2(M_1(k),P'_2) > x^*_2$.
Then $\gamma_2(M_1,P'_2) -c|P'_2|\leq \gamma_2(M_1,M_2)- c|M_2|$,
provided $c > \frac{3}{k^3k_2}$.
\end{lemma}

\begin{proof}
Let $L_0$ be the number of pure automata $(P_1^l)$ that $P_2$ fools.
Because $\gamma_2(M_1,P'_2) > x^*_2$ we have $L_0 \geq 1$.
If $P_2$ fools $P_1^l$, player 2's long-run average payoff is at most 4, the maximal payoff in the game.
If $P_2$ does not fool $P_1^l$, player 1's long-run average payoff is at most $x_2^*$.
The expected long-run average payoff of player 2 then satisfies
\[ \gamma_2(M_1,P'_2) \leq 4  \frac{L_0}{k_2} + x^*_2\frac{k_2-L_0}{k_2}
<  x_2^* + 3\frac{L_0}{k_2}. \]
By Corollary \ref{corollary5} we have $|P'_2| \geq L_0k^3$,
and therefore
\[ \gamma_2(M_1,P'_2) < x_2^* + 3\frac{L_0}{k_2} = x_2^* + 3\frac{L_0k^3}{k^3k_2}
\leq x_2^* + |P'_2| \times \frac{3}{k^3k_2}.\]
Therefore, as soon as $c > \frac{3}{k^3k_2}$ player 2 does not profit by this deviation.
\end{proof}

\bigskip

To summarize, given a feasible and an individually rational payoff vector $x^*$,
we choose $\eta \in (0,\min\{x_1^*-1,x_2^*-1\})$.
Let $c > 0$ be sufficiently small, and let $k=k_c$ satisfy
$\frac{3}{k^3k_2} < c < \frac{\eta}{3k^3}$.
Then the automata $(M_1(k),M_2(k))$ form a $c$-BCC equilibrium.
Since the size of the automata $M_1(k)$ and $M_2(k)$ are $k^3+2k^2+1$ and $k^3+2k^2+k+1$,
if for each $k \geq 1$ we set $\widehat c_k = \frac{4}{k^3k_2}$, then $\frac{3}{k^3k_2} < \widehat c_k < \frac{\eta}{2k^3}$ and
$\widehat c_k M_1(k)$ and $\widehat c_k M_2(k)$ are both smaller than $\frac{10}{k_2}$,
which goes to 0 as $k$ goes to infinity (and $\widehat c_k$ goes to 0).
It follows that $x^*$ is a BCC equilibrium payoff.

\section{The General Case}
\label{section general}

In Section \ref{section example} we proved Theorem \ref{th1} for the Prisoner's Dilemma.
In the present section we explain how the proof should be adapted to prove the result for arbitrary games.
In the play path $\omega^*$, the punishment phase, as well as the regular play,
are similar to those in Section \ref{section example}, and only the babbling phase significantly changes.

Assume w.l.o.g. that payoffs are bounded by 1, and let $x \in F\cap V$.
To rule out trivial cases, assume that each player has at least two actions.
The vector $x$ is a convex combination of {\em all} the entries in the payoff matrix
\[ \left| x - \sum_{a \in A} \alpha_a u(a)\right| \leq \ep, \]
where $(\alpha_a)_{a \in A}$ are non-negative numbers summing to 1.
In fact, by Caratheodory's Theorem, $x$ is a convex combination of three entries in the payoff matrix.
Instead of handling separately each of the alternative configurations of these three entries,
we find it simpler to handle the general case.

Fix $\ep > 0$, a natural number $k_0 > \frac{1}{\ep}$, and a natural number $k$.
Let $(k_a)_{a \in A}$ be a collection of positive integers such that
(a) $\sum_{a \in A} k_a = k_0$, and
(b) $\left| k_a - \alpha_a k_0 \right| \leq 1$.
Define the regular path
\[ \omega_0 = \sum_{a =(a_1,a_2)\in A} k_a \times (a_1,a_2). \]
Then the average payoff along $\omega_0$ is within $\ep$ of $x$.

For each $i=1,2$, denote by $l_i = |A_i|$ the number of actions of player $i$,
and by
$A_i = \{a_i^1,a_i^2,\ldots,a_i^{l_i}\}$ her actions.
Assume w.l.o.g. that $l_1 \leq l_2$,
and that $a_i^1$ is the min-max strategy of player $i$ against player $3-i$.

The play path $\omega^*$ is defined as follows:
\begin{eqnarray*}
\omega^* &=&
k^4 \times (a_1^1,a_2^1) + \left( \sum_{j=1}^{k^2} \sum_{m=1}^{l_1} k\times(a_1^m,a_2^m)\right)
+(k+1) \times (a_1^1,a_2^1)\\
&& + \left( \sum_{j=1}^{k} \sum_{m=l_1+1}^{l_2}k \times (a_1^1,a_2^m)\right) + (a_1^2,a_2^1) + \sum_{j=1}^\infty \omega_0.
\end{eqnarray*}
Both the punishment phase and the babbling phase are longer in this construction
than in the construction for the Prisoner's Dilemma,
yet the punishment phase is much longer,
to ensure that the payoff that results from a deviation is close to the min-max value in pure strategies.

The complexity of $\omega^*$ w.r.t. player 1 is $k^4 + k^3l_1 + 1$.
Indeed, as in Section \ref{section example}, the complexity of the prefix
$\omega^*(1) = k^4 \times (a_1^1,a_2^1) +
\left( \sum_{j=1}^{k^2} \sum_{m=1}^{l_1} k\times(a_1^m,a_2^m)\right)  + (k+1)\times(a_1^1,a_2^1)$
is $k^4 + k^3l_1 + 1$, and to implement the rest of the play $\omega^*$ player 1 can re-use states that were used to
implement $\omega^*(1)$.
One can verify that the complexity of $\omega^*$ w.r.t. player 2 is $k^4 + k^3l_1 + (k+1) + k^2(l_2-l_1)$.
A similar construction of automata $M_1$ and $M_2$ for the two players,
that re-uses states to implement the rest of $\omega^*$,
shows that $x$ is a BCC equilibrium payoff.

\end{document}